\documentclass{amsart}
\usepackage{epsf, amssymb, amsmath}

\newtheorem{prop}{Proposition}[section]
\newtheorem{theorem}[prop]{Theorem} 
\newtheorem{lemma}[prop]{Lemma}
\newtheorem{quest}{Question}

\newcommand{\R}{{\bf R}} 
\newcommand{\Rt}{{\bf R^{3}}}  
\newcommand{\Be}{{B_{\epsilon}}}
\newcommand{\epsfig}[1]{{\epsffile{#1.eps}}}

\title{Short Ropes and Long Knots}

\author{Jacob Mostovoy}

\address{Instituto de Matem\'{a}ticas (Unidad Cuernavaca),  
	Universidad Nacional Aut\'{o}noma de M\'{e}xico,  
	A.P. 273-3,  C.P. 62251, Cuernavaca, Morelos, MEXICO}  

\begin{document}
\maketitle

We study spaces of non-singular smooth embeddings of a closed interval into
$\R^3$. Our main result is a geometric interpretation of the Grothendieck group
of the monoid of knots.


\section{Ropes}

Let us begin with  definitions.
Fix two points $A$ and $B$ in $\Rt$. We shall always assume that
$A=(0,0,0)$ and $B=(1,0,0)$ so, in particular, the line $AB$ is the
$x$-axis and the length of the interval $[AB]$ is equal to $1$.
(We will write $[AB]$ for the closed interval between $A$ and $B$,
$(AB)$ for the corresponding open interval, and $AB$ for the line passing
through $A$ and $B$.)

A {\it rope} is a non-singular $C^1$-smooth embedding $r:[0,1]\to\Rt$
such that $r(0)=A$ and $r(1)=B$. Non-singularity here includes the
condition that the tangent vectors $\frac{dr}{dt}(0)$ and $\frac{dr}{dt}(1)$
are non-zero. 

The space of all ropes endowed with $C^1$-topology is denoted by $B_{\infty}$.
We will also consider, for any $\epsilon>0$, its subspaces $\Be$ that are 
formed by ropes whose length is strictly less than $1+\epsilon$. Each of 
these subspaces comes with a natural  basepoint --- the tight rope which is 
the embedding $t\to (t,0,0)$. Finally, we say that a rope is {\it short}, 
if its length is less than 3.
 
\begin{figure}[ht]
	\[\epsfig{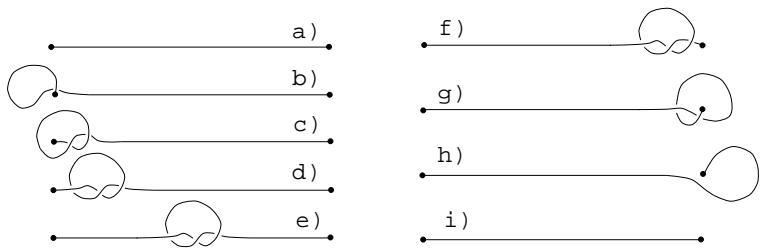}\]
	\caption{The image of a trefoil under 
					$b_{\epsilon}$.}\label{f:def-homo}
	\end{figure}

Even though ropes are knotty objects, spaces of ropes $\Be$ have only
one connected component. Indeed, any knotted rope can be ``undone'', 
i.e.\ deformed to the tight rope without increasing its length.
(Figure~\ref{f:def-homo} (f)-(i) illustrates how a knotted rope can be undone.)
However, the fundamental group of $\Be$ can be non-trivial and
we shall see that, for spaces of short ropes, it is closely related
to the semigroup of knots.

Unless stated otherwise, by ``knots'' we will mean ``long'' or 
``non-compact'' knots, i.e.\ nonsingular smooth embeddings $\R\to\Rt$
whose tangent vectors tend to $(1,0,0)$ at $\pm \infty$. These are
essentially knots in $S^3$ and their isotopy classes are in one-to-one
correspondence with the isotopy classes of ``round'' knots. Recall that the
isotopy classes of knots form a commutative monoid (i.e.\ a semigroup with a 
unit) $K$
under the connected sum. The monoid 
$K$ is freely generated by classes of prime knots, of which there are countably
many; the unknot, that is, the inclusion map of the $x$-axis into $\R^3$,
being the unit. We will often say ``knots''  for both geometric objects 
and their isotopy classes, it will always be clear from the context which we 
are talking about. 

To each $k\in K$ we can associate an element $b_{\epsilon}(k)\in\pi_1(\Be)$\
for any $\epsilon$ as follows. First by carrying the rope around the point $A$
we tie the knot $k$ on the rope near $A$. Then we push the knot
along the rope all the way to the point $B$ and throw it off the rope there.   
(See Figure~\ref{f:def-homo}.) This process defines a closed path in $\Be$
and as we will see later the homotopy class of this path corresponds to a 
well-defined element of $\pi_1(\Be)$. The precise definition of the map
$b_{\epsilon}$ will be given in Section~3.
It is easy to see that the map $b_{\epsilon}$ respects the connected sum of 
knots; that is to say that it is a homomorphism of the monoid $K$ to the 
group $\pi_1(\Be)$. 

Our main result is the following 
\begin{theorem}\label{t:main}
For any $0<\epsilon\leqslant 2$ the homomorphism $b_{\epsilon}:K\to\pi_1(\Be)$
is a group completion.
\end{theorem}
In other words, for any space of short ropes $\Be$ the group $\pi_1(\Be)$ is 
the Grothendieck group $\widehat{K}$ of the monoid $K$. It is easy to 
understand 
which loop corresponds to  $-k\in\widehat{K}$: it is defined in the same way
as $b_{\epsilon}(k)$, but we tie the knot at $B$ and push it to the left 
towards $A$.

\begin{figure}[ht]
	\[\epsfig{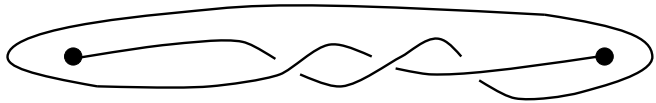}\]
	\caption{A rope which is not short.}\label{f:long-rope}
	\end{figure}

Interestingly, the situation changes dramatically as soon as $\epsilon$
becomes bigger than 2.
\begin{theorem}\label{t:long}
For any $\epsilon > 2$ the space $\Be$ is simply-connected.
\end{theorem}
The reason for such behaviour can be roughly explained as follows. To undo a  
knot which is tied on a short rope one needs to carry the rope either around
$A$ or $B$. However, if the rope is longer than 3, we can take the knot off 
the whole interval $[AB]$, see Figure~\ref{f:long-rope}.

At the moment we cannot say anything about the higher homotopy or homology 
groups of $\Be$ for any finite $\epsilon$. However, in the limit case 
$\epsilon=\infty$ an explicit deformation retraction of $B_{\infty}$ to the 
tight rope can be constructed. Thus we have the following:

\begin{theorem}\label{t:all}
The space of all ropes $B_{\infty}$ is contractible.
\end{theorem}

Our last result says that all spaces of short ropes have the 
same homotopy type:
\begin{theorem}\label{t:tight}
The natural inclusion  $\Be\hookrightarrow B_{\epsilon'}$
is a homotopy equivalence for all $0<\epsilon\leqslant\epsilon'\leqslant 2$.
\end{theorem}
This is all we know about ropes so far. This work is far from being a 
conclusive study and some questions are listed in the last section. 
The next section briefly 
explains where the idea of studying ropes comes from and in sections
3, 4, 5 and 6 we prove Theorems~\ref{t:main}, \ref{t:long}, \ref{t:all} and 
\ref{t:tight} respectively. 

\medskip

Much of the inspiration for this work came from discussions with
Sergei Chmutov, Sergei Duzhin, Elmer Rees and Simon Willerton who has also 
given some valuable comments on the preliminary version of this paper.


\section{Motivation: classifying spaces}

The motivation for studying ropes comes from the construction of a classifying
space for a topological monoid, described by M.C.McCord in \cite{McCord}.
Consider a space $BM$ whose points are configurations of particles on the
interval $[0,1]$ and all particles are labelled by non-zero elements of 
a monoid $M$. The topology is introduced in such a way that it agrees with the
topologies on [0,1] and the label space $M$; if two particles with coordinates
$x_1$ and $x_2$ (where $x_1<x_2$) and labels $h_1$ and $h_2$ respectively move
towards each other and eventually collide, they form a particle with the label
$h_{1}h_{2}$. The endpoints of the interval are ``sources'' of particles,
i.e.\ particles with arbitrary labels can appear from 0 and 1 and, conversely,
when a particle collides with 0 or 1 it disappears. The basepoint in $BM$
is chosen to be the empty configuration. (See \cite{McCord} for the 
precise definition.)   

\begin{theorem}\label{t:McCord}{\rm \cite{McCord}}
For $M$ a simplicial monoid, $BM$ is a classifying space of $M$. 
\end{theorem}
\noindent {\bf Remark.}
McCord's theorem as stated in \cite{McCord} deals with classifying spaces
of groups rather than monoids. McCord's construction, 
however, is exactly the same thing as Segal's construction of a 
classifying space for a category \cite{Seg} in the case when the category has 
only one object. 

\medskip

An example of this construction is the infinite symmetric product of a circle
$SP^{\infty}(S^1)$. Here particles are labelled by positive integers and labels
are just the multiplicities of points. The particles of $SP^{\infty}(S^1)$
live on a circle rather than on an interval; notice, however, that in McCord's 
construction points 0 and 1 can be identified. The infinite symmetric product
of $S^1$ is well-known to be topologically a circle, and this
is the classifying space for the monoid of non-negative integers.

If $G$ is a discrete group there is a map $\omega:G\to\Omega BG$ which 
induces an isomorphism $G\to\pi_1(BG)$.
The image of $g\in G$ under $\omega$ is the loop which at time $t$ is a 
configuration with a single particle which has coordinate $t$ and label $g$.
If $M$ is a discrete abelian monoid the map $\omega$ is also defined in the
same way and it induces the standard homomorphism of $M$ to its Grothendieck 
group $\widehat{M}$. The latter statement is a particular case of the Homology 
Group Completion Theorem (see \cite{BP} or \cite{MS}). Clearly, if 
$m_1,m_2\in M$, the element of $\pi_1(BM)$ which corresponds to the image of 
$m_1-m_2$ in $\widehat{M}$ can be represented by the loop 
$\omega(m_1)\overline{\omega}(m_2)$; here by $\overline{\omega}(m_2)$ we mean 
the loop $\omega(m_2)$ taken with the opposite parametrisation.

McCord's construction can, of course, be applied to the monoid $K$ of the 
isotopy classes of knots. The classifying space $BK$ can be then thought of
as an interval on which infinitesimally small knots are tied. Placing
the interval into $\Rt$ and replacing infinitesimally small knots by
knots of finite size we come to the notion of a rope; the map $b_{\epsilon}$
described in Section~1 is an analogue of the map $\omega$.
The condition $\epsilon\leqslant 2$ in this context means that knots that are
tied on a rope can be ``localised'' as particles. This is the main idea of 
the proof of Theorem~\ref{t:main}.


\section{Spaces of short ropes}
\begin{figure}[ht]
	\[\epsfig{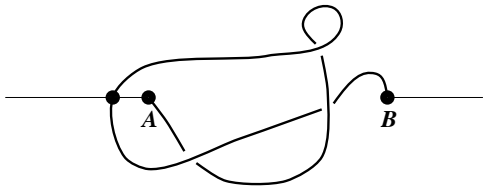}\]
	\caption{A knot extension of a rope.}\label{f:extension}
	\end{figure}

We define an {\it extension} of a rope $r$ as the map $\tilde{r}:\R\to\R^3$
which coincides with $r$ on the interval $[0,1]$ and such that
$\tilde{r}(t)=(t,0,0)$ for all $t\in (-\infty,0]\cup [1,\infty)$. 
An extension of a rope is a piecewise-smooth long knot which, however, 
can have points of self-intersection, see Figure~\ref{f:extension}. 
We say that an extension of a rope $r$ has no singularities  
to the left of $A$ if for all $t\in(-\infty,0]$ the point $\tilde{r}(t)$
is not a point of self-intersection. Similarly one defines what it means 
for an extension of a rope to have no singularities to the right of $B$.

Now we can give the precise definition of the map $b_{\epsilon}$.
An element of $\pi_1(\Be)$ can be represented by a loop on $\Be$, that is, by
a one-parameter family of ropes $L_T:[0,1]\to\Be$ such that 
$L_0=L_1={\rm tight\ rope}$. 
We say that a family  $L_T:[0,1]\to\Be$ is {\it generic} if there is a 
finite number of values of the parameter $T_i$ and $T'_j$ such that the 
ropes $L_{T_i}$ and $L_{T'_j}$ extend to knots with only one transversal 
double point to the left of $A$ or to the right of $B$ respectively, and 
for all other values of $T$ the rope $L_{T}$ extends to a genuine knot. 

The map $b_{\epsilon}$ assigns to a knot $k$ a one-parameter family of ropes
$L_T:[0,1]\to\Be$ which is generic in the above sense and such that:
\begin{itemize}
\item{there exists $X\in [0,1]$ for which $T_i<X<T_{j}^{'}$ for all $i,j$;}
\item{the knot extension of $L_X$ has isotopy class $k$;}
\item{$L_0=L_1={\rm tight\ rope}$.}
\end{itemize}   
\begin{lemma}\label{lemma:welldef}
The map $b_{\epsilon}:K\to\pi_1(\Be)$ is well-defined.
\end{lemma}
To prove this it is enough to show that if a loop $L_T$ on $\Be$
can be extended to the right of $B$ without singularities for any $T\in [0,1]$,
then it is contractible in $\Be$. This is an immediate corollary of
Lemma~\ref{l:WL-contraction} below.
\medskip

Assuming the truth of Lemma~\ref{lemma:welldef},  
we shall now prove that the map 
$\hat{b}_{\epsilon}$ of the Grothendieck group
$\widehat{K}$ into $\pi_1(\Be)$, which is induced by the map $b_{\epsilon}$
is an isomorphism for any $0<\epsilon\leqslant 2$.

\begin{prop}\label{prop:mono}
For $0<\epsilon\leqslant 2$ the map  
$\hat{b}_{\epsilon}:\widehat{K}\to\pi_1(\Be)$ is a monomorphism.
\end{prop} 
\begin{proof}
Let $p$ be the projection onto the line $AB$ and for a rope $r$ let 
$A(r)\subset AB$ be the subset of such points $t$ that the inverse image of 
$p\circ r(t)$ in $[0,1]$ consists of more than one point or that the tangent 
line to $r$ at $t$ is orthogonal to $AB$, see Figure~\ref{f:def} (a).
We say that a rope $r$ is {\it nice} if $A(r)$ is a union of finitely many
closed intervals and points. For example, analytic ropes are nice.

	\begin{figure}[ht]
	\[\epsfig{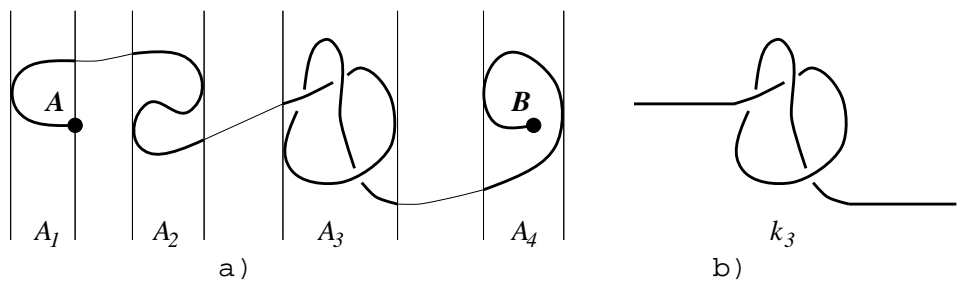}\]
	\caption{A nice rope and one of the corresponding knots.}\label{f:def}
	\end{figure}

Suppose $r$ is a nice rope and $A(r)=\bigcup A_i(r)$, where $A_i(r)$ are 
disjoint closed intervals or points (here $i$ belongs to some finite index 
set). To each $A_i(r)$ which does not contain $A$ or $B$ corresponds a knot 
$k_i(r)$ which is obtained by extending the segment of the rope $r$ 
which projects onto $A_i(r)$ to the
left and to the right
by rays parallel to the line $AB$ (see Figure~\ref{f:def} (b)). Clearly, if
$A_i(r)$ is a point, the corresponding $k_i(r)$ is a trivial knot. We say that
nice ropes $r_1$ and $r_2$ are of the same type if $A(r_1)$ can be identified
with $A(r_2)$ by some orientation-preserving self-homeomorphism of $AB$
which fixes $A$ and $B$, and if knots that correspond to $A_i(r_1)$ and 
$A_i(r_2)$ are the same.

Recall that a point in $BK$ is a collection of particles labelled by knots 
on a closed interval. Let us identify this interval with the interval $[AB]$.
 We say that a point $y\in BK$ is {\it subordinate}
to a nice rope $r$ if all particles of $y$ belong to $A(r)$ and if the 
sum of the coefficients of particles that are contained in $A_i(r)$ is exactly
$k_i(r)$. (If $A_i(r)$ contains $A$ or $B$ we allow the coefficients of particles
within $A_i(r)$ to add up to any knot.)

Let $S(r)\subset BK$ be subspace of points that are subordinate to the rope 
$r$. Clearly, if $r_1$ and $r_2$ have the same type, $S(r_1)$ and $S(r_2)$ 
are homeomorphic. Notice that $\epsilon\leqslant 2$ implies that the intervals
$A_i$ do not cover the whole $[AB]$ and, hence, for any short rope $r$ the 
space $S(r)$ is contractible.

Now suppose that a formal difference of knots $k_1-k_2\in\widehat{K}$
defines a contractible loop $\gamma_T$ in $\Be$. We can assume that $\gamma_T$
goes through nice ropes and that the standard loop 
$\omega_T(k_1)\overline{\omega_T}(k_2)$ which
represents $k_1-k_2$ in $BK$ (see Section 2) is subordinate to $\gamma_T$
at any value of the parameter $T$.
So there exists a map $F:[0,1]^2\to\Be$ which coincides with $\gamma_T$ on one 
side of the square and sends the rest of its boundary to the tight rope.
\begin{lemma}\label{l:fourier}
Without loss of generality we can assume that:
\begin{itemize}
\item{the image of $F$ consists of nice ropes only;}
\item{the square $[0,1]^2$ can be triangulated in such a way that the interiors
of all simplices of triangulation are mapped to ropes of the same type.} 
\end{itemize}
\end{lemma}  

The proof is purely technical. One can use Fourier expansions to 
replace the 2-dimensional family of ropes $F$ by an analytic 
family of analytic ropes; in this context the lemma is easily verified.
The details, in which neither God nor devil are to be found, are left to the 
reader.
 
Consider the subspace $S\subset [0,1]^2\times BK$ of all pairs $\{x,y\}$ where
$x\in [0,1]^2$ and $y\in S(F(x))$. Under the assumptions of 
Lemma~\ref{l:fourier} one can check that the projection $S\to [0,1]^2$
is a quasifibration with contractible fibres. (For the definition and 
properties of quasifibrations see \cite{DT}.) By the weak homotopy lifting
property of quasifibrations there is a map $F':[0,1]^2\to S$ which sends
$(T,0)\in [0,1]^2$ to $\{(T,0),\omega_T(k_1)\overline{\omega_T}(k_2)\}$ and 
any other point $x$ of the boundary to $\{x, *\}$; 
here $*\in BK$ is the basepoint. Combining $F'$ with the
projection $S\to BK$ we see that $\omega_T(k_1)\overline{\omega_T}(k_2)$ is 
contractible in $BK$ and this is possible only if $k_1=k_2$.
\end{proof}

\begin{prop}
For any $\epsilon >0$ the map  
$\hat{b}_{\epsilon}:\widehat{K}\to\pi_1(\Be)$ is an epimorphism.
\end{prop} 

\begin{proof}
Let $W_L\subset\Be$ be the
subspace of ropes which extend to knots without singularities to the left
of the point $A$ and let $W_R$ denote the subspace of ropes extending without 
singularities to the right of $B$.

\begin{lemma}\label{l:WL-contraction}
For any $\epsilon>0$ the subspaces $W_L$ and $W_R$ are contractible. 
\end{lemma} 

The proof of this lemma is rather technical and is better visualised than
verbalised. We postpone it till the end of the section.

\medskip
Ropes which lie in the intersection $W_L\cap W_R$ extend to non-singular
(apart from a possible discontinuity of the tangent vectors at $A$ and $B$)
knots. So $\pi_0(W_L\cap W_R)=K$, as it is easy to check that the restriction
on the length of the ropes does not influence the picture. As $W_L$ and
$W_R$ are contractible, the union $W_L\cup W_R$ has the homotopy type of the 
suspension on $W_L\cap W_R$. Consequently, $\pi_1(W_L\cup W_R)$ is the
free group, generated by the non-zero elements of $K$.

The complement of $W_L\cup W_R$ in $\Be$ has codimension 2, so the map
\[ \pi_1(W_L\cup W_R)\to\pi_1(\Be) \] 
induced by the inclusion $W_L\cup W_R\to\Be$
is onto. The problem now is to find the relations in $\pi_1(W_L\cup W_R)$.

Recall our definition of the map $b_{\epsilon}:K\to\pi_1(\Be)$. It is
clear from the definition that each knot defines, in fact, a loop in
$W_L\cup W_R$, as a rope in the process of deformation intersects the line
$AB$ first on the left of $A$ and then on the right of $B$. So there is
a well-defined map $K\to\pi_1(W_L\cup W_R)$ which, however, is not a 
homomorphism. One can easily identify the image of this map: 
a knot is mapped to the corresponding generator of the group 
$\pi_1(W_L\cup W_R)$.

We know that the composite map 
\[ K\to\pi_1(W_L\cup W_R)\to\pi_1(\Be)\]
is a homomorphism. In case $\epsilon\leqslant 2$ this immediately 
determines the set of
relations we are looking for: the generators must commute as $K$ is abelian,
and $k_1k_2$ (the product in $\pi_1(W_L\cup W_R)$) must be equal to the 
connected
sum $k_1\# k_2$. These relations define the Grothendieck group $\widehat{K}$
and, as we have seen above, for $\epsilon\leqslant 2$ the group 
$\pi_1(\Be)$ cannot be smaller.
\end{proof}

In case  $\epsilon > 2$ there may be some additional relations. And indeed,
in the next section we shall see that in this case $\pi_1(\Be)=0$.

\begin{proof}[{\em Proof of Lemma~\ref{l:WL-contraction}.}]
The spaces $W_L$ and $W_R$ are homeomorphic so it suffices to
verify the statement of the lemma for $W_L$.

The deformation retraction of $W_L$ to the tight rope is done in two stages.
Let $W_L^0\subset W_L$ be the subset of ropes whose tangent vector at $A$ 
is $(a,0,0)$ for some $a>0$; the extensions of all ropes in $W_L^0$ to the 
left can be parametrised so as to be smooth. First we will show how to deform
$W_L$ into $W_L^0$ and then prove that $W_L^0$ can be deformed to the
tight rope.   

Choose the spherical coordinates $(r,\phi,\theta)$ in $\Rt$ with the centre at
$A$ and $\theta=0$ being the positive half of the $x$-axis in  $\Rt$.
Consider the family of maps
\[d_T: (r,\phi,\theta)\to (r,\phi,\theta\cdot (1-\frac{5T}{6}))\]
with $T\in [0,1]$. 
Away from the non-positve part of the 
$x$-axis and for each $T\in [0,1]$ the map $d_T$ is a diffeomorphism onto
its image. Moreover, $d_T$ does not increase lengths and thus we have the 
corresponding deformation of the space of ropes
\[D_T: W_L\times [0,1]\to W_L.\]

The effect of the map $D_1$ is that all ropes in $W_L$ are pushed into the 
cone $\theta<\pi/6$.
In order to deform them further to $W^0_L$ we need to ``squeeze'' the tip of 
this cone, see Figure~\ref{f:stage1}. The subtle point here is that is has to 
be done while keeping the lengths of ropes under control.

 	\begin{figure}[ht]
 	\[\epsfig{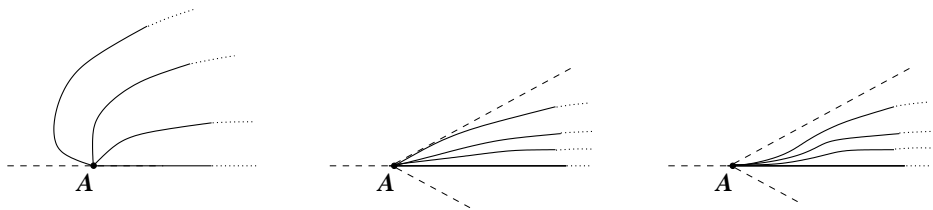}\]
 	\caption{The deformation of $W_L$ to $W_L^0$ near the 
 	point $A$.}\label{f:stage1}
 	\end{figure}

Define the distance $\rho(r_1,r_2)$ between two ropes $r_1$ and $r_2$ by
\[ \rho(r_1,r_2) = 
\max_{t\in [0,1]} (|r_1(t)-r_2(t)|+|\frac{d}{dt}(r_1(t)-r_2(t))|). \]
Together with the Euclidean metric on $\R$ the distance function 
$\rho$ gives rise to a distance function on 
$B_{\epsilon}\times\R$. Identifying the line $AB$ with $\R$ in the 
obvious way, that is, $A=0$ and $B=1$, we obtain a distance function
on  $\Be\times AB$.
 
Let $E\subset\Be\times AB$ be the subspace of pairs $(r,x)$ such that
$(p\circ r)^{-1}(x)$ consists of only one point at which the angle between
the tangent vector to $r$ and the line $AB$ is less than $\pi/4$. Notice
that for any rope $r\in D_1(W_L)$ the pair $(r,0)$ lies in $E$.

Define $\delta_1(r)$ to be the infimum of the distance between $(r,0)$
and the complement of $E$ in $\Be\times AB$, 
and let  $\delta_2(r)$ be equal to
$1+\epsilon-l(r)$, where $l(r)$ is the length of the rope $r$. On the subspace
$D_1(W_L)$ both $\delta_1(r)$ and $\delta_2(r)$ are continuous and 
strictly 
positive functions of $r$, so if we set
\[\delta(r)=\frac{1}{5}\min (\delta_1(r),\delta_2(r))\]
the function $\delta(r)$ is also positive and continuous on  $D_1(W_L)$.

Let us fix a smooth function $f: [0,+\infty)\to\R$ such that
$f(0)=0$, $f(x)=1$ for $x\geqslant 1$ and $0<\frac{df}{dx}<2$ for
$x<1$. For $r\in D_1(W_L)$ and $T\in [0,1]$ consider the function
$f_{r,T}(x)=Tf(\frac{x}{\delta(r)})+(1-T)$. Notice that $f_{r,T}(x)$ is 
continuous in $r$ as $\delta(r)$ is positive and continuous on $D_1(W_L)$. 
Then the deformation retraction 
\[D'_T: D_1(W_L)\times [0,1]\to D_1(W_L) \]
can be defined on a rope $r=(r_x,r_y,r_z)$ as
\[D'_T(r) = 
(r_x, f_{r,T}(r_x)r_y, f_{r,T}(r_x)r_z).\]

The boundedness of $\frac{df}{dx}$ and the
condition that $f(0)=0$ imply that at $T=1$ all ropes are carried into
ropes whose tangent vectors point in the direction of the $x$-axis.
It remains to check that $D'_T(r)$ does not increase the lengths of ropes
too much.

Notice that under $D'_T$ a rope $r$ changes only near the point $A$, namely
within the region $0<x<\delta(r)$. Recall that by definition we have
$\delta(r)<\delta_1(r)$ so for all $0<x<\delta(r)$ the pair $(r,x)$ belongs to
$E$. This, in particular, means that for $0<x<\delta(r)$ the rope $r$ can
be re-parametrised as $(x,r_y(x),r_z(x))$ with 
$(\frac{d}{dx}r_y(x))^2+(\frac{d}{dx}r_z(x))^2<1$. It follows that the part of
$r$ within the region  $0<x<\delta(r)$ has length between $\delta(r)$ and 
$\delta(r)\sqrt{2}$.

Now,
\[\left\lvert\frac{d}{dx}(f_{r,T}r_y)\right\rvert\leqslant 
\left\lvert\frac{d}{dx}(f_{r,T})\right\rvert |r_y|+
 |f_{r,T}|\left\lvert\frac{d}{dx}r_y\right\rvert
\leqslant \left\lvert\frac{2T}{\delta(r)}\right\rvert |\delta(r)|+
1\cdot 1\leqslant 3\]
and the same inequality holds for $|\frac{d}{dx}(f_{r,T}r_z)|$.
Hence, the length of the part of $D'_T(r)$ 
contained in the  region $0<x<\delta(r)$ is
bounded by $\delta(r)\sqrt{1+3^2+3^2}< 5\delta(r)$. It follows that $D'_T$
can only increase the length of a rope $r$ by less than 
$5\delta(r)\leqslant\delta_2(r)= (1+\epsilon)-l(r)$ which means that 
the total length of $D'_T(r)$ is less than $1+\epsilon$ for all $T$ and $r$.

Let us now describe the second stage, which is the deformation 
$D''_T$ of $W_L^0$ to the 
tight rope. For each $r$ we construct $D''_T(r)$ 
(with $D''_1$ being the identity map 
and $D''_0$ --- the map to the tight rope) 
in several steps. Without loss of generality we will assume that all ropes 
are parametrised by length (up to a constant factor).

First we cut the rope $r$ at the value of the parameter $T$, i.e.\ consider
the embedding $r:[0,T]\to\Rt$. If the projection $p(r(T))$ of the point
$r(T)$ onto $AB$ lies to the right of the point $B$, we ``squeeze the rope'' 
linearly by a map $s:\Rt\to\Rt$ given by
\[ (x,y,z)\to(\frac{x}{|p(r(T))|}, y,z)\]
so that $p(s(r(T)))=B$.

If $p(r(T))$ lies to the left of $B$ we ``move the rope to the right'' by a
translation $s'$, such that  $p(s'(r(T)))=B$, see Figure~\ref{f:move}.
(The interval between $A$ and $(1-p(r(T)),0,0)$ is filled in with a segment
of a straight 
line.) One can check that when $p(r(T))\to B$ both $s$ and $s'$ tend to the
identity map.

 	\begin{figure}[ht]
 	\[\epsfig{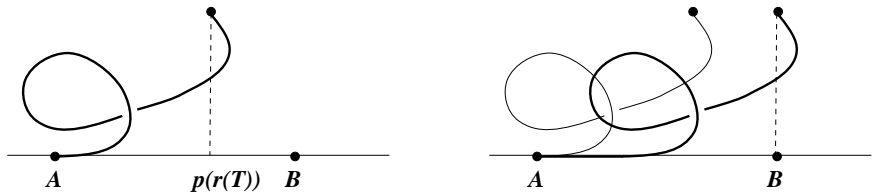}\]
 	\caption{``Moving a rope to the right''.}\label{f:move}
 	\end{figure}

One of the endpoints of the resulting curve is $A$; let us denote the
other endpoint by $B'$. Then
the second step is a shift in the planes orthogonal to $AB$ given by 
\[ (x,y,z)\to (x,y,z)-x^2\cdot (B'-(1,0,0)).\]
After this shift the curve obtained is a rope $r'$ (or, more precisely, 
an image of a rope) as the ends of it coincide
with $A$ and $B$ respectively. The length of $r'$, however, may exceed
$1+\epsilon$ so we need the third step. It is a transformation induced by 
squeezing $\Rt$ to the
line $AB$ with the help of a continuous function $h_{r,T}$ such that 
$0< h_{r,T}\leqslant 1$:
\[ (x,y,z)\to (x, yh_{r,T}, zh_{r,T}).\]

The function  $h_{r,T}$ is chosen to be equal to 1 if the length of $r'$ is
not greater than the length of $r$; otherwise we define $h_{r,T}$ by the
condition that the rope $(r'_x,r'_y h_{r,T},r'_z h_{r,T})$ has the same length
as $r$.

Finally, we parametrise the rope by
length (up to a constant factor) and the desired deformation $D''_T$ is the 
composition of all the above transformations.

Now it can be checked directly that the map $D''_T:W_L^0\times [0,1]\to W_L^0$ 
is continuous and that $D''_0(W_L^0)$ is the tight rope.
\end{proof}


\section{Spaces of long ropes}

In the previous section we have seen that for any $\epsilon>0$ the map
$\widehat{K}\to \pi_1(\Be)$ is onto. Here we will show that if $\epsilon>2$
this map is a Vassiliev invariant of order 1. It is a well-known fact that
all Vassiliev knot invariants of order 1 
are constants; so this will imply that $\pi_1(\Be)=0$ for $\epsilon>2$.

 	\begin{figure}[ht]
 	\[\epsfig{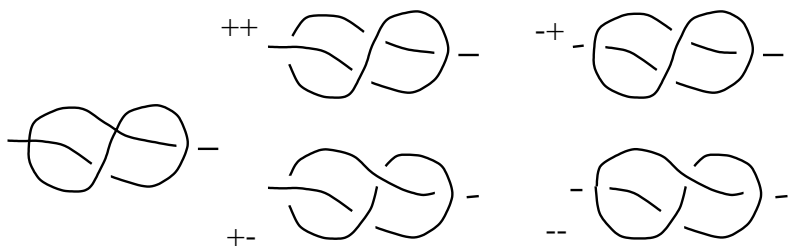}\]
 	\caption{A singular knot $k$ 
 	and knots $k_{++}, k_{+-},k_{-+}$ and $k_{--}$.}
 	\label{f:resolution}
 	\end{figure}

\begin{figure}[ht]
 	\[\epsfig{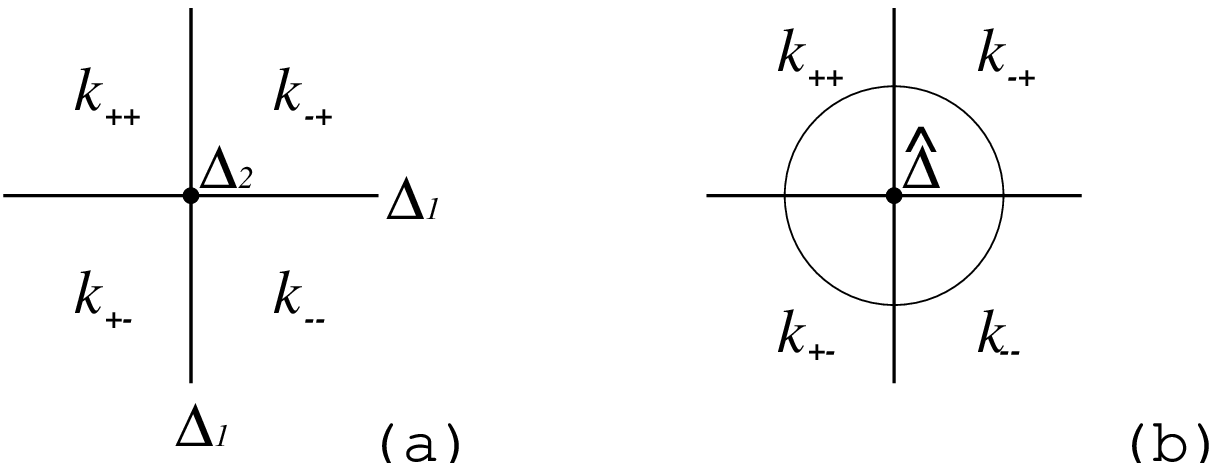}\]
 	\caption{}\label{f:delta}
 	\end{figure}

For basic facts about Vassiliev invariants we refer the reader to \cite{BL}, 
\cite{Stan} or \cite{Vas}. 
Here we recall the definition of a Vassiliev invariant of order 1. 
In the space of
all smooth maps $\R\to\Rt$ that coincide with a chosen line outside some 
finite interval there is a discriminant $\Delta$ which is formed by 
non-embeddings. The space of knots is then the complement of $\Delta$ in the
space of all maps. The stratum $\Delta_1$ of $\Delta$ which has codimension 1 
is formed
by knots with generic double points and the stratum $\Delta_2$ of 
codimension 2 is formed
by knots with 2 double points. A neighbourhood of a point on $\Delta_2$
is pictured on Figure~\ref{f:delta} (a). Here $k_{++}$, $k_{+-}$,  
$k_{--}$ and $k_{-+}$ are the knots obtained by resolving the singularities
of a knot with two double points. An example of such resolution is given 
on Figure~\ref{f:resolution}. A function 
\[ v_1:K\to G \] which takes values in some abelian group $G$ 
is a Vassiliev invariant of order 1 if for all points of $\Delta_2$
\[ v_1(k_{++})-v_1(k_{+-})+v_1(k_{--})-v_1(k_{-+})=0.\]
The same definition is valid for round knots, i.e.\ embeddings $S^1\to\Rt$.

Any knot invariant can be extended from $K$ to $\widehat{K}$ by linearity.
In particular, a Vassiliev invariant of order 1 on $\widehat{K}$ is a
linear extension of a Vassiliev knot invariant of the same order. It is
clear that a knot invariant is identically zero if and only if its extension
to $\widehat{K}$ is.

We have a very similar picture in the space of ropes. Indeed, the complement 
$\widehat{\Delta}$ of $W_L\cup W_R$ in $\Be$ has codimension 2 and a 
neighbourhood of a generic point of $\widehat{\Delta}$ is pictured on 
Figure~\ref{f:delta} (b). Here the 
vertical line is formed by ropes which extend to a knot with one double point
to the left of $A$ and the horizontal line is formed by ropes which extend with
a double point to the right of $B$. The ropes that extend to genuine knots
are labelled by the type of their knot extensions: $k_{++}$,  $k_{+-}$,  
$k_{--}$ and $k_{-+}$.

A circle around $\widehat{\Delta}$ represents the zero element in $H_1(\Be)$
(which is equal to $\pi_1(\Be)$ as $\pi_1(\Be)$ is abelian), however it 
might not represent the zero element in $H_1(W_L\cup W_R)$. Recall that 
$W_L\cup W_R$ is a suspension on the space of knots (that are subject to some
length restriction). So $H_1(W_L\cup W_R)$ is a free abelian group,
generated by non-zero elements of $K$. Identifying the knots $k_{\pm,\pm}$
with the corresponding generators of $H_1(W_L\cup W_R)$ we see that a circle 
around $\widehat{\Delta}$ defines, up to sign, the element 
$k_{++}-k_{+-}+k_{--}-k_{-+}$ in $H_1(W_L\cup W_R)$; and so 
$k_{++}-k_{+-}+k_{--}-k_{-+}\in\widehat{K}$ is sent to 0 in $\pi_1(\Be)$.

It remains to show that if the knots $k_{++}$,$k_{+-}$,$k_{--}$ and $k_{-+}$ 
are found as knot types in the neighbourhood of some point of the codimension-2
stratum in $\Delta$, they appear in the same cyclic order as types of
knot extensions near some generic point of $\widehat{\Delta}$ when 
$\epsilon>2$.

\begin{figure}[ht]
 	\[\epsfig{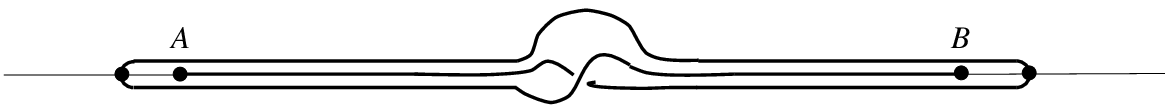}\]
 	\caption{}\label{f:twodblepts}
 	\end{figure}

Here, for once, we will make use of round knots. Take a round knot with 2
double points $x_1$ and $x_2$ and let $k_{++}$, $k_{+-}$, $k_{--}$ and 
$k_{-+}$ be the resolutions of its singularities. Consider a segment 
of the knot that connects
$x_1$ and $x_2$ and choose a point on this segment. If we take this point 
to be the infinity in $S^3$, we get a long knot $f$ such that there are no
double points to the left of the smaller value of the parameter $t_1$ that 
corresponds to $x_1$ and to the right of the larger value of the parameter
$t_2$ that corresponds to $x_2$. Let $a$ be slightly larger than $t_1$
and $b$ slightly smaller that $t_2$. We can deform our knot in such a way that
$f(a)=A$, $f(b)=B$ and that it is an extension of a rope of length less than
$1+\epsilon$ for any $\epsilon>2$, see Figure~\ref{f:twodblepts}. Clearly, 
this  provides us with a rope in $\widehat{\Delta}$ in whose neighbourhood
the knot extensions have types $k_{++}$,$k_{+-}$,$k_{--}$ and $k_{-+}$. 
This proves Theorem~\ref{t:long}. 

\section{The space of all ropes}

Here we construct an explicit deformation retraction 
of the space $B_{\infty}$ onto the tight rope.  
The deformation $\delta_T:B_{\infty}\times [0,1]\to B_{\infty}$ is the 
composition of the following transformations:

Similarly to the proof of Theorem~\ref{t:all}, first we cut the rope at the 
value of the parameter $T$, i.e.\ consider the embedding $r:[0,T]\to\Rt$.
The next step is a homothety with the centre at $A$:
\[ (x,y,z)\to\frac{(x,y,z)}{|r(T)|}.\]
After this we rotate $\Rt$ around $A$ so that the ``free end of the rope'',
i.e.\ the point $\frac{r(T)}{|r(T)|}$ is moved to $B$. The rotation $R(r,T)$ 
is determined from the condition that the derivative $\frac{dR(r,t)}{dt}(T)$ 
is an infinitesimal rotation around the axis $r(T)\times\frac{dr}{dt}(T)$ of 
magnitude $|r(T)\times\frac{dr}{dt}(T)|\cdot|r(T)|^{-2}$ and such that 
$R(r,0)={\bf Id}$.

It is a straightforward check that $\delta_T$ is a continuous deformation
retraction  with $\delta_1$ being the identity map and 
$\delta_0$ --- the map to the tight rope.


\section{Tightening the ropes}

The proof of Theorem~\ref{t:tight} resembles in spirit the proofs of 
Lemma~\ref{l:WL-contraction} and Theorem~\ref{t:all}: we construct an 
explicit deformation retraction $B_{\epsilon '}\to\Be$ for any 
$0<\epsilon<\epsilon'\leqslant 2$.
In other words, we will show how to tighten all ropes in $B_{\epsilon '}$
simultaneously. 

First of all let us introduce some notation: $l(r)$ will stand for the length 
of the rope $r=(r_{x}(t),r_{y}(t),r_{z}(t))$, and $l_{x}(r)$ and $l_{yz}(r)$ 
are the lengths of paths $(r_{x}(t),0,0)$ and $(0,r_{y}(t),r_{z}(t))$ 
respectively.
 
The retraction consists of two steps. First we reduce $l_{yz}(r)$ by squeezing
$\Rt$ to the line $AB$ with the help of some continuous function $h_{r,T}$
with the arguments $r\in B_{\epsilon '}$ and $T\in [0,1]$:
\[ (x,y,z)\to (x, yh_{r,T}, zh_{r,T}), \]
(compare this with the proof of Lemma~\ref{l:WL-contraction}). The second step
reduces $l_{x}(r)$; the deformation of the space of ropes in this case is 
also induced by a family of deformations of $\Rt$. Here they are of the form  
\[ (x,y,z)\to (x\phi_{r,T}(x), y, z), \]
where $\phi_{r,T}(x)$ is a family of $C^1$-smooth monotonic functions 
$\R\to\R$
which depends continuously on parameters $r$ and $T$.

The function $h_{r,T}$ is chosen as follows. Let 
$\psi_{\tau}:B_{\epsilon'}\times [0,1)\to B_{\epsilon'}$ be the deformation
given by
\[ \{r,\tau\}\to (r_{x}, (1-\tau)r_{y}, (1-\tau)r_{z}).\]
It is clear that for any $r\in B_{\epsilon'}$, apart from the tight rope, 
$l_{yz}(\psi_{\tau}(r))$ and $l(\psi_{\tau}(r))$ are 
decreasing functions of $\tau$ and for all ropes $l_{x}(\psi_{\tau}(r))$ 
does not depend on $\tau$. 
Notice that the function $l_{yz}(\psi_{\tau}(r))$ is linear in $\tau$ 
and continuous in $r$, and $\lim_{\tau\to 1}{l_{yz}(\psi_{\tau}(r))}=0$ 
for any $r\in B_{\epsilon'}$. So for all ropes apart from the tight rope 
we can define $f_1(r)$ as the minimal value of $\tau$ such that 
$l_{yz}(\psi_{\tau}(r))\leqslant \frac{\epsilon}{2}$ and 
$f_2(r)$ as the minimal 
value of $\tau$ such that 
$l_{yz}(\psi_{\tau}(r))+l_{x}(\psi_{\tau}(r))\leqslant l(r)$. 
Now set $f(r)=1$ if $r$ is the tight rope and
$f(r)=\max\{f_1(r),f_2(r)\}$ otherwise. The continuity of $f_1(r)$ and $f_2(r)$
away from the tight rope follows from the linearity 
of $l_{yz}(\psi_{\tau}(r))$. 
It can also be checked directly that as $r$ tends to the tight rope, 
$f_1(r)$ tends to 0 and $f_2(r)$ tends to 1. (In fact, $f_1$ is identically 
zero in some neighbourhood of the tight rope.)
Thus, $f(r)$ is also a continuous function 
$B_{\epsilon'}\to [0,1]$ and $f(r)=1$ if and only if $r$ is the tight rope.
Finally, we define $h_{r,T}$ as
\[ h_{r,T}= 1 - Tf(r).\]

If $H:B_{\epsilon'}\to B_{\epsilon'}$ denotes the map
\[ r\to (r_{x},  r_{y}h_{r,1}, r_{z}h_{r,1})\]
it is clear from the above construction that 
$l_{yz}(H(r))\leqslant\frac{\epsilon}{2}$ and $l_{yz}(H(r))+l_{x}(H(r))\leqslant l(r)$
for all $r$.

The second step is more involved. 

Recall that in the proof of the Propostion~\ref{prop:mono} we defined for 
each rope $r$ a subset $A(r)\in AB$.
Let $Z(r)$ be the complement of $A(r)\cap [AB]$ in $[AB]$.
In other words, $Z(r)$ is the subset of interval $[AB]$ formed by such points
$x$ that $(p\circ r)^{-1}(x)$ consists of only one point at which 
$\frac{dr_{x}}{dt}$ is not zero. 
The key point in what follows is
that for a short rope $Z(r)$ is a non-empty open subset of $[0,1]$
so it has non-zero Lebesgue measure. 

The total length of the projection of $r$ onto $AB$
can be written as a sum $l_{x}(r)=l_{A}(r)+l_{Z}(r),$ where $l_{A}(r)$
and $l_{Z}(r)$ are lengths of the parts of $p(r)$ that lie in $A(r)$ and
$Z(r)$ respectively. Clearly, $l_{Z}(r)$ is just the Lebesgue measure of
$Z(r)$.

Suppose that we have a family of $C^1$-smooth monotonic functions 
$\phi_{r,T}:\R\to\R$
which depends continuously on parameters $r\in B_{\epsilon'}$ and $T\in [0,1]$
and satisfies the following conditions:

(a) $\phi_{r,T}(0)=0$ and $\phi_{r,T}(1)=1$ for all $r$ and $T$;

(b) $\frac{d}{dx}\phi_{r,T}(x)<1$ for all $x\in A(r)$ and any $T\in (0,1]$;

(c) $\frac{d}{dx}\phi_{r,1}(x)\leqslant\frac{\epsilon}{4l_{A}(r)}$ 
for any $x\in A(r)$.\\
Then we can define a homotopy 
\[ \Phi_{T}: H(B_{\epsilon'})\times [0,1] \to  B_{\epsilon'} \]
which is given by 
\[ \{r(t),T\}\to (\phi_{r,T}(r_{x}(t)),r_{y}(t),r_{z}(t)). \]
\begin{lemma}
The homotopy $\Phi_{T}$ is well-defined and $\Phi_{T}(H(\Be))\subset\Be$
for any $T$.
\end{lemma}
\begin{proof}
The condition (a) above means that $\Phi_{T}$ takes ropes to ropes; so
in order to show that  $\Phi_{T}$ is well-defined we need to check that it 
does not increase the lengths of ropes too much.

It is clear that the Lebesgue measure of $A(r)\cap [AB]$ is equal to
$(1-l_{Z}(r))$. 
Let $n(x)$ be the number of inverse images of $p\circ r$ at $x\in AB$.
The function $n(x)$ is finite almost everywhere on $AB$ and 
$n(x)-1$ is nonnegative on $A(r)$. Integrating $n(x)$ over $A(r)$ we obtain
$l_{A}(r)$, and the sum of integrals, 
\[ I(r)=\int_{A(r)\cap[AB]}(n(x)-1)dx+
\int_{A(r)\backslash(A(r)\cap[AB])}n(x)dx, \]
is equal to $l_{A}(r)-1+l_{Z}(r)$.
The condition (b) above implies that
\[I(\Phi_{T}(r))\leqslant I(r)\] 
so 
\[ l_{x}(\Phi_{T}(r))=l_{A}(\Phi_{T}(r))+l_{Z}(\Phi_{T}(r))
\leqslant l_{A}(r)+l_{Z}(r)=l_{x}(r).\]
Now recall that for any $r\in B_{\epsilon'}$ the sum of lengths of the 
projections $l_{x}(H(r))+l_{yz}(H(r))$ is less than or equal to $l(r)$.
In particular, for all $r\in H(B_{\epsilon'})$ we have
$l_{x}(r)+l_{yz}(r)<1+\epsilon'$. 
It is easy to see  $\Phi_{T}$ does not change $l_{yz}(r),$
so the length of $\Phi_{T}(r)$ is bounded by $l_{x}(r)+l_{yz}(r)<1+\epsilon'$
and this means that $\Phi_{T}$ is well-defined.

Similarly, the inequality $l_{x}(r)+l_{yz}(r)<1+\epsilon$ holds
for all $r\in H(B_{\epsilon})$  and for all such ropes
$l(\Phi_{T}(r))$ is bounded by $1+\epsilon$. This means that 
$\Phi_{T}(H(\Be))\subset\Be$.
\end{proof}

\begin{lemma}
The map $\Phi_{1}\circ H: B_{\epsilon'}\to \Be$ 
is a homotopy equivalence. 
\end{lemma}
\begin{proof}
Recall that $l_{yz}(r)\leqslant\frac{\epsilon}{2}$ for any 
$r\in H(B_{\epsilon'})$. The condition (c) implies that
\[ l_{A}(\Phi_{1}(r))\leqslant l_{A}(r)\cdot\frac{\epsilon}{4l_{A}(r)}= 
\frac{\epsilon}{4},\]
hence the length of the rope $\Phi_{1}(r)$ is bounded by
\[ l_{A}(\Phi_{1}(r))+l_{Z}(\Phi_{1}(r))+l_{yz}(\Phi_{1}(r))\leqslant
\frac{\epsilon}{4}+1+\frac{\epsilon}{2}<1+\epsilon.\]
So $\Phi_T$ deforms $H(B_{\epsilon'})$ into $\Be$. Moreover, 
for any $T$ the map $\Phi_T$ takes $H(B_{\epsilon})$ into $\Be$ and thus
$\Phi_{1}\circ H$ is a homotopy equivalence.
\end{proof}

To finish the proof we need to construct the family of functions $\phi_{r,T}$.

Let ${\mathfrak A}\subset B_{\epsilon'}\times AB$ be the subset
\[ {\mathfrak A}=\{ r,u |\ r\in B_{\epsilon'}, u\in AB, u\in A(r).\]
Recall that in the proof of Lemma~\ref{l:WL-contraction} we have constructed 
a distance function on the spaces $B_{\epsilon}\times AB$.
For each rope $r$ define a  function $g_{r}(x)$ to be the infimum of the 
distance from $(r,x)$ to the subset ${\mathfrak A}$. 
It is clear that the family 
of functions $g_{r}(x)$ is continuous both in $x$ and $r$.
Notice that $g_{r}(x)=0$ for any $x\in A(r)$ and $g_{r}(x)>0$ for any 
$x\in Z(r)$. This implies, in particular, that for any short rope
$r$ the integral $\int_{0}^{1}g_{r}(y)dy$ is greater than zero.

Let ${\tilde{\phi}}_{r,T}(x)$  be equal to 
\[ \frac{\int_{0}^{x}1+Tg_{r}(y)dy}{\int_{0}^{1}1+Tg_{r}(y)dy}. \]
For any $x\in A(r)$
$$ \frac{d}{dx}{\tilde{\phi}}_{r,T}=
\left({\int_{0}^{1}1+Tg_{r}(y)dy}\right)^{-1} < 1,
\leqno{(\ast)}$$
so if we set 
\[ \beta(r)= \max\left( 0, \left[\frac{4l_{A}(r)}{\epsilon}-1\right]\cdot
\left[{\int_{0}^{1}g_{r}(y)dy}\right]^{-1}\right) \] 
and ${\phi}_{r,T}={\tilde{\phi}}_{r,T\beta(r)}$ the conditions (a), (b) and 
(c) are satisfied.

Indeed, (a) follows straight from the definitions and (b) is a direct 
consequence of $(\ast)$ above. For $x\in A(r)$ we have
\[ \frac{d}{dx}\phi_{r,1}(x)=\frac{d}{dx}\tilde{\phi}_{r,\beta}(x)=
\left(1+\beta(r){\int_{0}^{1}g_{r}(y)dy}\right)^{-1}\leqslant 
\frac{\epsilon}{4l_{A}(r)} \]
by $(\ast)$ and the definition of $\beta(r)$. 
This verifies (c) and all that is left is to check the continuity of the 
family ${\phi}_{r,T}$. This immediately follows from the continuity of.
${\tilde{\phi}}_{r,T}(x)$ and $\beta(r)$ in all variables.

\section{Final remarks and questions}

It is certainly interesting that  spaces of short ropes provide us with
a geometric interpretation of the Grothendieck group of the monoid of
knots. However, a stronger question can be asked. 

Recall that behind the spaces of short ropes there is a construction of a 
classifying space. Long knots form an $H$-space and it is not too hard 
to modify this $H$-space to make it associative, i.e.\ to make it a 
monoid.  
\begin{quest}
Is the space $\Be$ a classifying space for long knots when $\epsilon\leqslant 2$?
\end{quest} 
For a simplicial monoid $M$ the space $\Omega BM$ of loops
on the classifying space of $M$ can be considered as a ``homological
group completion'' of $M$, see \cite{BP} or \cite{MS}. So there is a chance
that the topology of the spaces of short ropes is related to the topology of
the space of knots in a rather direct way.

Another question concerns knot invariants. Having got a generic 
loop $L_T$ in a space of short ropes we can find out what element in 
$\widehat{K}$ it corresponds to as follows. 

Suppose that the  knot extensions of $L_T$ have double points to the left of 
$A$ at the values of the parameter $T_i$, and let $\delta_i k\in\widehat{K}$ 
be the formal differences of the types of knot extensions near $T_i$. Then the
loop $L_T$ represents the formal difference 
$\sum_{T_i}\delta_i k\in\widehat{K}$.
(This easily follows from the fact that $W_L\cup W_R$ is a suspension on the 
space of non-singular knot extensions.) Similarly one can calculate the value
of any additive knot invariant on any element of $\pi_1(\Be)$. 

This, however, is not interesting in the sense that we learn nothing new 
about knots.
\begin{quest}
Which knot invariants can be defined geometrically on the level of ropes?  

\end{quest} 

Vassiliev invariants have shown up in our constructions, though in a rather
silly way. It would be good to understand if there is a deeper connection.
It is interesting that with the help of Vassiliev invariants on can construct  
groups of knots in a less simple-minded way than taking the Grothendieck group.
These are Gusarov's groups of $n$-equivalence classes of knots, see \cite{Gus}
 and \cite{NgSt}. Also, cobordism classes of knots form a group, \cite{FoxMil}.
Clearly, there are homomorphisms from $\widehat{K}$ to all of these groups;
it is not clear, however, if these homomorphisms can be interpreted 
geometrically via ropes.


\end{document}